\pdfoutput=1
\documentclass[12pt,a4paper]{article}
\usepackage{amsmath,amssymb,amsthm,eucal,MnSymbol,amscd,relsize}

\usepackage[text={453pt,686pt},centering]{geometry}	%
\usepackage[noconfig]{refstyle}				%
\usepackage{microtype,amsfonts}				%
\usepackage[usenames,dvipsnames,svgnames,table]{xcolor}	%
\usepackage[numbers,sort&compress]{natbib}		%
\usepackage[british]{babel}					%
\usepackage{hyphenat}					%
\allowdisplaybreaks[4]					%
\usepackage[colorlinks, citecolor=blue, linkcolor=blue, urlcolor=Maroon, filecolor=Maroon, linktocpage=true]{hyperref}

\font\tenrm=cmr10

\font\bigss=cmssdc10 scaled 2300

\font\cmsslll=cmss10 at 14 pt  
\renewcommand{\a}{\alpha}  
\renewcommand{\b}{\beta}

\newcommand{\g}{\gamma}

\renewcommand{\l}{\lambda}

\newcommand{\s}{\sigma}

\newcommand{\G}{\Gamma}

\newcommand{\bR}{\mathbb{R}}  
\newcommand{\bZ}{\mathbb{Z}}

\renewcommand{\gg}{\mathfrak{g}}

\newcommand\Sp{\mathrm{Sp}}  
\newcommand\GL{\mathrm{GL}}  
  
\newcommand\SO{\mathrm{SO}}  
\newcommand\SU{\mathrm{SU}}  
\newcommand\U{\mathrm{U}}

\newcommand{\p}{\partial}

\DeclareMathOperator\tr{tr}  
  
\DeclareMathOperator\End{End}  
\DeclareMathOperator\vol{vol}  
  
\DeclareMathOperator\ad{ad}

\newtheorem{Th}{Theorem}%
\newtheorem{Prop}[Th]{Proposition}  
\newtheorem{Cor}[Th]{Corollary}  
\newtheorem{Lem}[Th]{Lemma}  

\theoremstyle{definition} %
\newtheorem{Def}[Th]{Definition}  
\newtheorem{Ex}[Th]{Example}

\theoremstyle{remark}

\newtheorem*{remark}{Remark}

\newcommand{\bt}{\begin{Th}\ \ }  
\newcommand{\et}{\end{Th}}  
\newcommand{\bp}{\begin{Prop}\ \ }  
\newcommand{\ep}{\end{Prop}}  
\newcommand{\bc}{\begin{Cor}\ \ }  
\newcommand{\ec}{\end{Cor}}  
\newcommand{\bl}{\begin{Lem}\ \ }  
\newcommand{\el}{\end{Lem}}  
\newcommand{\bd}{\begin{Def}\ \ }  
\newcommand{\ed}{\end{Def}}  
\newcommand{\bex}{\begin{Ex}\ \ }  
\newcommand{\eex}{\end{Ex}}  

\newcommand{\pf}{\begin{proof}}
\newcommand{\epf}{\end{proof}}

\newcommand{\be}{\begin{equation}}  
\newcommand{\ee}{\end{equation}}

\newcommand{\arr}{\begin{array}{rlll}}  
\newcommand{\ea}{\end{array}}  
\newcommand{\bea}{\begin{eqnarray}}  
\newcommand{\eea}{\end{eqnarray}}  
\newcommand{\bean}{\begin{eqnarray*}}  
\newcommand{\eean}{\end{eqnarray*}}  
\let\eqref=\relax
\newref{eq}{name={eq.~},Name={Eq.~},names={eqs.~},Names={Eqs.~},rngtxt={-},refcmd=(\ref{#1})}
\newref{f}{name={footnote~},Name={Footnote~},names={footnotes~},Names={Footnotes~}}
\newref{app}{name={appendix~},Name={Appendix~},names={appendixes~},Names={Appendixes~}}
\newref{tab}{name={table~},Name={Table~},names={tables~},Names={Tables~}}
\newref{ch}{name={chapter~},Name={Chapter~},names={chapters~},Names={Chapters~}}
\newref{sec}{name={section~},Name={Section~},names={sections~},Names={Sections~}}
\newref{fig}{name={figure~},Name={Figure~},names={figures~},Names={Figures~}}
\newref{thm}{name={Theorem~},Name={Theorem~},names={Theorems~},Names={Theorems~}}
\newref{prop}{name={Proposition~},Name={Proposition~},names={Propositions~},Names={Propositions~}}
\numberwithin{equation}{section}
\newcommand{\RR}{\mathbb{R}}
\newcommand{\CC}{\mathbb{C}}
\newcommand{\PP}{\mathbb{P}}
\newcommand{\ZZ}{\mathbb{Z}}

\newcommand{\V}{\mathcal{V}}

\DeclareMathOperator\grad{grad}
\DeclareMathOperator\Isom{Isom}
\def\a{\alpha}
\def\b{\beta}

\def\1{{\bar 1}}
\def\2{{\bar 2}}
\def\3{{\bar 3}}
\def\4{{\bar 4}}

\def\and{\quad\textrm{and}\quad}

\newcommand\varpm{\mathbin{\vcenter{\hbox{%
  \oalign{\hfil$\scriptstyle+$\hfil\cr
          \noalign{\kern-.5ex}
          $\scriptscriptstyle({-})$\cr}%
}}}}

\newcommand{\resumetocwriting}{%
  \addtocontents{toc}{\protect\setcounter{tocdepth}{\arabic{tocdepth}}}}
\begin{document}  
\hypersetup{pageanchor=false}
\begin{titlepage}

\begin{flushright}
\normalsize{\texttt{ZMP-HH/17-13}}\\
\normalsize{\texttt{HBM 653}}\\
\end{flushright}

\rightline
{June 27, 2018}

\vskip 2.0 true cm  
\begin{center}  
{\bigss  Left-invariant Einstein metrics on $\mathlarger{\mathlarger{\mathlarger{\boldsymbol{S^3 \times S^3}}}}$
}

\vskip 1.5 true cm   
{\cmsslll Florin\ Belgun$^\S$, Vicente\ Cort\'es$^\P$, Alexander\ S.\ Haupt$^\P$, and David\ Lindemann$^\P$} \\[2em] 
{\tenrm $^\S$``Simion Stoilow" Institute of Mathematics of the Romanian Academy\\
Calea Grivitei 21,
Sector 1, 010702 Bucharest, Romania\\
\texttt{florin.belgun@uni-hamburg.de}}\\[2em]  
{\tenrm $^\P$
Department of Mathematics and Center for Mathematical Physics\\
University of Hamburg, Bundesstr. 55, 
D-20146 Hamburg, Germany\\  
\texttt{\{vicente.cortes, alexander.haupt, david.lindemann\}@uni-hamburg.de}
}\\[1em]  

\vspace{2ex}

\end{center}

\vspace*{1em}
\baselineskip=18pt  
\begin{abstract}  
\noindent
The classification of homogeneous compact Einstein manifolds in dimension six is an open problem. 
We consider the remaining open case, namely left-invariant Einstein metrics $g$ on $G = \mathrm{SU}(2) \times \mathrm{SU}(2) = S^3 \times S^3$. 
Einstein metrics are critical points of the total scalar curvature functional for fixed volume.
The scalar curvature $S$ of a left-invariant metric $g$ is constant and can be expressed as a rational function in the parameters determining the metric.
The critical points of $S$, subject to the volume constraint, are given by the zero locus of a system of polynomials in the parameters.
In general, however, the determination of the zero locus is apparently out of reach. Instead, we consider the case where the isotropy group $K$ of $g$ 
in the group of motions is non-trivial.
When $K\not\cong \mathbb{Z}_2$ we prove that the Einstein metrics on $G$ are given by (up to homothety) either the standard metric or the nearly K\"ahler metric, based on representation-theoretic arguments and computer algebra.
For the remaining case $K\cong \mathbb{Z}_2$ we present partial results.

\end{abstract}

\end{titlepage}
\hypersetup{pageanchor=true}

\setcounter{tocdepth}{2} %

\newpage 

\section{Introduction and motivation}

In this paper we continue the study of homogeneous compact Einstein manifolds in six dimensions, see~\cite{NikRod2,NikRod} and references therein. 
The main progress is that we are able to treat the case where the stabilizer is finite rather than continuous.
Recall that, originating from the theory of general relativity, an Einstein manifold is defined to be a (pseudo-)Riemannian\footnote{In this work, we only consider the Riemannian case.} manifold $(M,g)$ whose Ricci tensor $Ric_g$ satisfies 
\be\eqlabel{Einstcond}
  Ric_g = \l g \; ,
\ee
for some constant $\l\in\RR$ called Einstein constant. The trace of this equation yields
\be\eqlabel{Einstcondtrace}
 S = n \l \; ,
\ee
where $S$ denotes the scalar curvature of $g$ and $n := \dim M$.

In~\cite{NikRod2,NikRod} a partial classification of such manifolds was obtained, stating that a simply connected six-dimensional homogeneous compact Einstein manifold is either a symmetric space or isometric, up to
multiplication of the metric $g$ by a constant, to one of the following manifolds:
(1) $\CC\PP^3 = \frac{\Sp(2)}{\Sp(1)\times \U(1)}$ with the squashed metric,
(2) the Wallach space $\SU(3)/T_{\text{max}}$ with the standard metric or with the K\"ahler metric, or
(3) the Lie group $\SU(2) \times \SU(2) = S^3 \times S^3$ with some left-invariant Einstein metric. 
Here and in the following we will consider $S^3$ as the group of unit quaternions. 
Hence, in order to complete the classification it is necessary to classify left-invariant Einstein metrics on $S^3 \times S^3$ (up to isometry).
The latter classification problem is still open. However, progress can be achieved by assuming additional symmetries of the metric $g$ (see, for instance, \thmref{NikonorovThm,Z2xZ2} below).

For left-invariant Einstein metrics on $S^3 \times S^3 =: G$, up to changing the metric by an isometric left-invariant metric, we have that~\cite[corollary on page 23]{DAtriZiller}
\be \eqlabel{Isom0}
  L_G \subset \Isom_0 (G,g) \subset L_G \cdot R_G \cong (G \times G)/\{(z,z) \mid z\in Z(G) \} \; ,
\ee
where $\Isom_0 (G,g)$ is the connected isometry group of some left-invariant metric $g$ on $G$, $L_G$ ($R_G$) is the group of left (right) translations and
$Z(G) \cong \bZ_2 \times \bZ_2$ denotes the center of $G$.
The right-hand side of \eqref{Isom0} contains the group of inner automorphisms 
\be
  \mathrm{Inn} (G)=C_G := \{ C_a \mid a \in G \} \subset L_G \cdot R_G \; , 
\ee
where $C_a$ denotes conjugation by $a$, that is 
\be
  C_a : G \to G \; , \qquad x \mapsto a x a^{-1} \; .
\ee
Hence, the isotropy group of the neutral element $e \in G$ in $\Isom_0 (G,g)$ is given by
\be\eqlabel{K_0def}
  \Isom_0 (G,g) \cap C_G =: K_0 \; ,
\ee
which is the maximal connected subgroup of the Lie group
\be\eqlabel{Kdef}
  \Isom (G,g) \cap C_G =: K \; .
\ee
In~\cite{NikRod}, a classification was achieved\footnote{For the sake of accurateness, we note that the last equation system on page 377 of~\cite{NikRod} contains a minuscule typo, which has however no influence on other parts of the presentation. Namely, in the third line the third term from the left should read $tuv(w-ut)$ instead of $tv(w-ut)$. } for the case that $K$ (or, equivalently, $K_0$) contains a $U(1)$ subgroup. This is summarized in the following theorem.
\begin{Th}[Nikonorov-Rodionov~\cite{NikRod}] \thmlabel{NikonorovThm} 
Let $g$ be a left-invariant Einstein metric on $G:=S^3 \times S^3$. If $K$, as defined in~\eqref*{Kdef}, contains a $U(1)$ subgroup, 
then $(G,g)$ is homothetic to $(G,g_{can})$ or $(G,g_{NK})$, where $g_{can}$ and $g_{NK}$ are the standard metric and the nearly K\"ahler metric, respectively.
\end{Th}
The two metrics $g_{can}$ and $g_{NK}$ are the only known Einstein metrics on $S^3\times S^3$ up to isometry and scale. It is also known
that these metrics are rigid. This follows from \cite[Proposition 4.8]{K} and \cite[Theorem 5.1]{MS}, respectively for the product metric and the nearly K\"ahler metric.
It is worth noting that $g_{can}$ is also right-invariant and, thus, invariant under the full adjoint group $\mathrm{Ad} (G) = \SO (3) \times \SO (3)$. 
The nearly K\"ahler (or Jensen's~\cite{Jensen}) metric $g_{NK}$ is only invariant under
the image $\mathrm{Ad} (\SU(2)_d) = \SO (3)_d :=\{ (a,a) \mid a\in \SO(3) \} \subset \SO (3) \times \SO(3)$ of the diagonal $\SU(2)$-subgroup 
$\SU (2)_d := \{ (a,a)\in G \mid a\in \SU (2)\}$ under the adjoint representation $\mathrm{Ad}=\mathrm{Ad}^G$ of $G$.

\Thmref{NikonorovThm} covers the case $\dim K \geq 1$. To complete the classification it remains to consider the case $\dim K = 0$, that is, the case where $K$ is a finite group. 
This is equivalent to requiring that $\Isom_0 (G,g)=G$, in which case the group of motions (that is orientation preserving
isometries) is given by 
\be \Isom^+ (G,g) = K  \ltimes G,\ee
where $K$ is a finite group of inner automorphisms of $G$. 
Analyzing this case is the goal of the present paper. Our main results can be summarized as follows.
\begin{Th}\thmlabel{Z2xZ2} Let $g$ be a left-invariant Einstein metric on $G$ that is invariant under a non-trivial finite subgroup $\G \subset \mathrm{Ad} (G)$ such that
$\G\not\cong \bZ_2$. Then $(G,g)$ is homothetic to $(G,g_{can})$ or $(G,g_{NK})$.
\end{Th}
The proof of this theorem requires a case-by-case analysis and concludes in \secref{Z2timesZ2}.
The case $\G = \bZ_2$ is considerably more complicated to analyze and in addition qualitatively novel features arise in the intermediate steps of the calculation. As a consequence, only partial results are available at this point. The Einstein condition on a left-invariant Riemannian metric on $G$ that is invariant under $\G = \bZ_2$ leads to a system of 12 coupled polynomial equations of degree 6 in 12 unknowns (see \secref{Z2}). 

Albeit solving the system is apparently out of reach with current technology, it is possible to analyze the space of solutions.
Whereas the systems of polynomial equations solved in the course of the proof of \thmref{Z2xZ2} 
have only a finite number of solutions, passing from groups $\G$ of order $\ge 2$ to $\G=\bZ_2$ leads to infinitely many
solutions.
\bp\proplabel{Z2contfam}
The system~\eqref*{Z2eqsys} of polynomial equations that describes
left-invariant Einstein metrics on $G$ invariant under a subgroup $\bZ_2 \subset \mathrm{Ad} (G)$ has 
continuous families of (real) solutions. 
\ep
However, all solutions of the system~\eqref*{Z2eqsys} which we have found so far are homothetic to $g_{can}$ or to $g_{NK}$, as we will 
explain now.
We have analyzed in more detail the aforementioned system of polynomial equations by holding fixed the value of the Lagrange multiplier\footnote{In our conventions, the Lagrange multiplier $\mu$ is related to the Einstein constant $\l$ via $\mu = -2 \l$ (see \secref{gen}).} $\mu$ of the variational problem (see \secref{gen}). Indeed, fixing $\mu$ eliminates it from the system, which can then be fully solved for the remaining variables. Of particular interest are the values $\mu=-1$ and $\mu=-5/(3 \sqrt{3})$ corresponding to the known solutions $(G,g_{can})$ and $(G,g_{NK})$, respectively. For these two values of $\mu$ we obtain continuous families of solutions, for other values of $\mu$ there are no solutions known. Furthermore we obtain the following result.
\bp\proplabel{Z2lambdafixed}
Let $g$ be a left-invariant Riemannian metric on $G$ that is invariant under a subgroup $\bZ_2 \subset \mathrm{Ad} (G)$. For $\mu=-1$, all solutions to the variational problem~\eqref*{var3}, which is equivalent to $g$ being Einstein with Einstein constant $\l = - \mu/2$, are isometric to a multiple of $g_{can}$. For $\mu=-5/(3 \sqrt{3})$, all solutions to the variational problem~\eqref*{var3} are isometric to a multiple of $g_{NK}$.
\ep

We end the introduction with some remarks, highlighting the relevance of six-dimen\-sional Einstein manifolds in the context of high energy physics.
Compact six-dimensional Einstein manifolds, in particular homogeneous spaces, feature prominently in various physical applications located mostly in the realm of string theory and its low-energy limit supergravity, as explained below.

Firstly, Einstein manifolds play a role in the AdS/CFT correspondence (see, for example,~\cite{Biquard:2005yq} and references therein). The conjecture asserts that string-/M-theory backgrounds of the form $\mathrm{AdS}_d \times M$, where $\mathrm{AdS}_d$ denotes $d$-dimensional anti-de Sitter space and $M$ is a compact Einstein manifold, should have an associated dual description as a conformal field theory on the $(d-1)$-dimensional boundary of $\mathrm{AdS}_d$. For example, type IIA superstring theory on $\mathrm{AdS}_4 \times \mathbb{CP}^3$ plays a role in the AdS$_4$/CFT$_3$ duality~\cite{ABJM}. Besides $\mathbb{CP}^3$ also $S^2 \times \mathbb{CP}^2$, $S^2 \times S^4$, and $S^2 \times S^2 \times S^2$ feature as possible compact Einstein six-manifolds in Freund-Rubin compactifications~\cite{FreundRubin} of type IIA supergravity to $\mathrm{AdS}_4$~\cite{Romans}. In the case of massive type IIA supergravity there are Freund-Rubin backgrounds of the form $\mathrm{AdS}_4 \times M$, where $M$ can be either $\mathbb{CP}^3$, the six-sphere $S^6$, the Grassmann manifold $SO(5)/(SO(2) \times SO(3))$, or one of the product spaces $S^3 \times S^3$, $S^2 \times \mathbb{CP}^2$, $S^2 \times S^4$, or $S^2 \times S^2 \times S^2$~\cite{Romans}.

Secondly, in (warped) flux compactifications of ten-dimensional string theory to four dimensions, the requirement of unbroken residual supersymmetry of the low-energy effective theory forces the internal six-dimensional manifold to admit an $\SU(3)$-structure~\cite{Strominger:1986uh, LopesCardoso:2002vpf, Becker:2003yv, Becker:2003gq, LopesCardoso:2003dvb, LopesCardoso:2003sp, Becker:2003sh, Becker:2003dz, Gray:2012md} (for reviews on the subject, see also, for example,~\cite{Grana:2005jc,Wecht:2007wu,Douglas:2006es,Blumenhagen:2006ci,Samtleben:2008pe}). Of particular interest are the cases where the $\SU(3)$-structure is nearly K\"ahler~\cite{Chatzistavrakidis:2008ii, Chatzistavrakidis:2009mh, Manousselis:2005xa, Lechtenfeld:2010dr, Klaput:2011mz, Chatzistavrakidis:2012qb, Gemmer:2012pp} or half-flat~\cite{Gurrieri:2002wz, Gurrieri:2004dt, Micu:2004tz, Gurrieri:2005af}. 
The (strict) nearly K\"ahler condition implies that the underlying Riemannian six-manifold is Einstein.
Besides the nearly K\"ahler metric, the product metric is an example of a left-invariant Einstein metric compatible with a left-invariant half-flat $\SU(3)$-structure~\cite{SchulteHengesbach:PhDThesis, SchulteHengesbach:2010jg}, see also~\cite{Madsen:2012by}. It is an open problem whether these are the only examples of such metrics on $S^3 \times S^3$ up to homothety.
For compactifications of heterotic supergravity with first-order $\a'$-corrections included, particular types of higher-dimensional Yang-Mills instantons arise as additional ingredients in the compactification set-up~\cite{Strominger:1990et,Harvey:1990eg,Khuri:1993ii,Gunaydin:1995ku,Loginov:2008tn,Harland:2011zs,Gemmer:2012pp,Klaput:2012vv,Haupt:2014ufa}. Consequently, instanton solutions of this type have been constructed, for example on cylinders, cones, and sine-cones over homogeneous compact nearly K\"ahler six-manifolds~\cite{Harland:2009yu, Gemmer:2011cp, Charbonneau:2015coa}.
The last subject is related to the topic of Hitchin flows over manifolds with half-flat $\SU(3)$-structure and other $G$-structures~\cite{Hitchin2001,Cortes2009}.
Note also that inhomogeneous compact nearly K\"ahler six-manifolds have recently been described in~\cite{Cortes2014} (locally homogeneous examples) and~\cite{Foscolo2015} (cohomogeneity one examples).

\vspace{\baselineskip}
\noindent\textit{Acknowledgements.} We thank Jos\'e V\'asquez for initial collaboration and discussions 
during later stages of the work.  We also thank Klaus Kr\"oncke and Yuri Nikonorov for helpful discussions. 
We gratefully acknowledge IT support from the IT-Group at the Department of Mathematics of the University of Hamburg, the HPC-Team at the RRZ of the University of Hamburg, and the Magma group of the University of Sydney.
This work was supported by the German Science Foundation (DFG) under the Collaborative Research Center (SFB) 676 ``Particles, Strings and the Early Universe''.

\section{Preliminaries}\seclabel{gen}

Finding Einstein metrics, that is finding solutions of \eqref{Einstcond}, can be reformulated as a variational problem~\cite{NikRod,Jensen,Besse,WangZiller}. Namely, a Riemannian metric $g$ on a compact orientable manifold $M$ solves \eqref{Einstcond} if and only if it is a critical point of the total scalar curvature functional, also known as the Einstein-Hilbert functional,
\be
 S_{\text{EH}} [g] = \int_M S \, \vol_g \; ,
\ee
subject to the volume constraint $\V := \int_M \vol_g = \V_0$, where $\V_0$ is a positive constant. Here, $\vol_g$ is the metric volume form on $(M,g)$. 

The volume constraint can be incorporated into the variational procedure by means of the method of Lagrange multipliers. Instead of directly varying $S_{\text{EH}} [g]$, we consider variations of
\be
 \tilde{S}_{\text{EH}} [g, \nu] = S_{\text{EH}} [g] - \nu (\V - \V_0) \; ,
\ee
where $\nu$ is a Lagrange multiplier. The vanishing of the variation of $\tilde{S}_{\text{EH}} [g, \nu]$ with respect to $g$ and $\nu$ yields
\be\eqlabel{var1}
 \grad_g S_{\text{EH}} [g] = \nu \grad_g \V \and \V = \V_0 \; ,
\ee
respectively. Here, $\grad_g$ denotes the variation with respect to the metric $g$. Plugging in the definitions of $S_{\text{EH}} [g]$ and $\V$, we obtain from the first equation in~\eqref*{var1}
\be
 \frac{S}{2} g - Ric_g = \frac{\nu}{2} g \; .
\ee
Comparing this to the Einstein condition~\eqref*{Einstcond} and using \eqref{Einstcondtrace} determines the Einstein constant $\l$ in terms of $\nu$, namely
\be
 \l = \frac{\nu}{n-2} \; ,
\ee
for $n>2$.

When $(M,g) = (G,g)$ is a compact Lie group $G$ (or more generally a unimodular Lie group, see~\cite[Theorem 1]{Jensen}) with left-invariant Riemannian metric $g$, simplifications occur in the general considerations above. In particular the scalar curvature $S$ is constant. Hence, $S_{\text{EH}} [g] = S\, \V$ and \eqref{var1} becomes
\be\eqlabel{var2}
 \grad_g S = -\frac{2 \l}{\V} \grad_g \V \and \V = \V_0 \; ,
\ee
which is equivalent to the Einstein condition~\eqref*{Einstcond} for metrics of unit volume.

Notice that a left-invariant Riemannian metric $g$ on $G$ is equivalent to a scalar product on the Lie algebra $\gg$ of $G$, which, for simplicity, we denote again by $g$.
Further specializing to $G = S^3 \times S^3$ and following~\cite{NikRod}, we consider the Lie algebra $\gg = \mathfrak{su}(2) \oplus \mathfrak{su}(2)$ with scalar product $Q(\cdot, \cdot) = -1/2 B(\cdot, \cdot)$, where $B(X,Y)=\tr(\ad(X)\ad(Y))$ is the Killing form of $\gg$. Any other scalar product $g$ can be obtained from $Q$ via $g(\cdot, \cdot) = Q(L \cdot, \cdot)$ for some ($Q$-symmetric) positive definite  endomorphism $L\in\End(\gg)$. In this way, the space of left-invariant Riemannian metrics is parameterized by the space 
\be 
 P(\gg) := \{ L\in\End(\gg) \mid \text{$L$ positive definite} \} \; .
\ee
Starting from $Q$ and some $Q$-orthonormal basis $(\mathbf{E}, \mathbf{F})$ of $\gg$, where $\mathbf{E} := (E_1, E_2, E_3)$, $\mathbf{F} := (F_1, F_2, F_3)$ are oriented orthonormal bases of the two $\mathfrak{su}(2)$-factors, we parameterize the space $P(\gg)$ by considering a change of basis from $(\mathbf{E}, \mathbf{F})$ to some $g$-orthonormal basis $(\mathbf{X}, \mathbf{Y})$ via
\be\eqlabel{genbasischange}
  (\mathbf{X} , \mathbf{Y} ) = (\mathbf{E} , \mathbf{F} ) A^T \; , \qquad A \in \GL (6,\RR) \; .
\ee
The matrix $A$ describing the change of basis satisfies $A^T A = L^{-1}$. We can choose $(\mathbf{X} , \mathbf{Y} )$ such that $A$ can be represented as~\cite{NikRod}
\be\eqlabel{genbasischangematrix}
 A = \begin{pmatrix} D & 0 \\ W & \tilde{D} \end{pmatrix} \; , \qquad\text{where}\qquad 
 D = \begin{pmatrix} a & 0 & 0 \\ 0 & b & 0 \\ 0 & 0 & c \end{pmatrix} \; , \quad
 \tilde{D} = \begin{pmatrix} d & 0 & 0 \\ 0 & e & 0 \\ 0 & 0 & f \end{pmatrix} \; , \quad
 W = \begin{pmatrix} x & u & v \\ \a & y & w \\ \b & \g & z \end{pmatrix} \; ,
\ee
such that $a, \ldots, f$ are positive parameters, whereas the components of $W$ are arbitrary real parameters.

Henceforth, we choose $\V_0 = \int_{G} \vol_Q = \left(2^3\mathrm{Vol}_{\overline{g}_{can}}(S^3)\right)^2=2^8\pi^4$,
%\begin{equation*}
%	\V_0 = \int_{G} \vol_Q = \left(2^3\mathrm{Vol}_{\overline{g}_{can}}(S^3)\right)^2=2^8\pi^4,
%\end{equation*}
where $\overline{g}_{can}$ denotes the canonical metric on $S^3\subset\mathbb{R}^4$. Note that $g_{can}= \overline{g}_{can}\oplus\overline{g}_{can}$, which together with $G=4g_{can}$ and $\mathrm{Vol}_{4\overline{g}_{can}}(S^3)=2^3\mathrm{Vol}_{\overline{g}_{can}}(S^3)$ explains the formula for $\mathcal{V}_0$.

The scalar curvature $S$ and the volume $\V=V\, \V_0$ of $g$ can be expressed as polynomials in the parameters $(a, \ldots, f, x, y, z, u, v, w, \a, \b, \g)$, namely
\begin{align}
 S &= a^2 +b^2 +c^2 + d^2 + e^2 + f^2 + x^2 + y^2 + z^2 +u^2 +v^2 +w^2 + \a^2 + \b^2 + \g^2 \nonumber \\
  &\quad - \frac12 \left\{ \vphantom{\left[ \left(\frac{de}{f}\right)^2 \right]} a^2 b^2 c^{-2} + b^2 c^2 a^{-2} + c^2 a^2 b^{-2} + d^2 e^2 f^{-2} + e^2 f^2 d^{-2} + f^2 d^2 e^{-2} \right. \nonumber \\ 
  &\quad + \left( \frac{a^2}{c^2} + \frac{c^2}{a^2} \right) (u^2 + y^2 + \g^2) + \left( \frac{a^2}{b^2} + \frac{b^2}{a^2} \right) (v^2 + w^2 + z^2) + \left( \frac{b^2}{c^2} + \frac{c^2}{b^2} \right) (x^2 + \a^2 + \b^2) 
  \nonumber \\
  &\quad + a^{-2} \left[ \left( uw - vy - \frac{de}{f} \b \right)^2 + \left( v\g - uz - \frac{df}{e} \a \right)^2 + \left( yz - w\g - \frac{ef}{d} x \right)^2 \right]
  \nonumber \\
  &\quad + b^{-2} \left[ \left( v\a - xw - \frac{de}{f} \g \right)^2 + \left( xz - v\b - \frac{df}{e} y \right)^2 + \left( w\b - z\a - \frac{ef}{d} u \right)^2 \right]
  \nonumber \\
  &\quad + c^{-2} \left. \left[ \left( xy - u\a - \frac{de}{f} z \right)^2 + \left( u\b - x\g - \frac{df}{e} w \right)^2 + \left( \a\g - y\b - \frac{ef}{d} v \right)^2 \right] \right\} \eqlabel{Sgen}
\end{align}
and
\be
 V\V_0 = (\det A)^{-1}\V_0 = (abcdef)^{-1}\V_0 \; ,
\ee
respectively~\cite{NikRod}. 
Einstein metrics then correspond to critical points of $S$ given by~\eqref*{Sgen} subject to the volume constraint $V=(abcdef)^{-1}=1$, that is, to solutions of
\be\eqlabel{var3}
 \nabla S = \mu \nabla V \and V = (abcdef)^{-1} = 1 \; ,
\ee
where $\nabla$ is the standard gradient in the parameter space $\left(\mathbb{R}_{>0}\right)^6\times \mathbb{R}^9 \subset \bR^{15}$ with the coordinates $(a, \ldots, f, x, y, z, u, v, w, \a, \b, \g)$ and $\mu$ is a Lagrange multiplier.

\begin{remark}
The relation between the Lagrange multiplier $\mu$ and the Einstein constant $\l$ can be clarified by comparing~\eqref*{var3} with~\eqref*{var2}. The first equation in~\eqref*{var3} can be written as
\be
 \frac{\p S}{\p A_{ij}} = \mu \frac{\p V}{\p A_{ij}} \; ,
\ee
where $A_{ij} = (A)_{ij}$ are the components of the matrix $A$. Using $L = (A^T A)^{-1}$ and the chain rule, we obtain
\be\eqlabel{remark_mu_lambda1}
 \sum_{k,l} \frac{\p L_{kl}}{\p A_{ij}} \frac{\p}{\p L_{kl}} (S - \mu V) = 0 \; .
\ee
For a $Q$-orthonormal basis, $g=L$ and hence the first equation in~\eqref*{var2} becomes
\be
 \frac{\p S}{\p L_{ij}} = -\frac{2 \l}{\V_0} \frac{\p \V}{\p L_{ij}} \; .
\ee
Inserting this into \eqref{remark_mu_lambda1} and using $\V=V\, \V_0 = \sqrt{\det L} \, \V_0$, we find
\be
 (\mu+2 \l) \tr\left( \frac{\p L}{\p A} A^T A \right) = 0 \; .
\ee
For the first factor inside the trace, we compute $\frac{\p L}{\p A} = - (A^T A)^{-1} (A + A^T) (A^T A)^{-1}$. Hence, $\tr\left( \frac{\p L}{\p A} A^T A \right) = - \tr ( A^{-1} + (A^T)^{-1})$ and after evaluating the trace using~\eqref*{genbasischangematrix}, we finally arrive at
\be
 (\mu+2 \l) \left( \frac1a + \frac1b + \frac1c + \frac1d + \frac1e + \frac1f \right) = 0 \; .
\ee
Since $a, \ldots, f$ are positive parameters, we conclude that $\mu=-2 \l$.
\end{remark}

We end this section by observing that the expression for $S$ as given in~\eqref*{Sgen} can be cast into a simpler form. This can be achieved by means of the following coordinate transformation $(a, \ldots, f, x, y, z, u, v, w, \a, \b, \g) \to (A, \ldots, F, X, Y, Z, U, V, W, \mathcal{A}, \mathcal{B}, \mathcal{C})$ of $\left(\mathbb{R}_{>0}\right)^6\times \mathbb{R}^9$,
\begin{align}
a&=\sqrt{BC}\; , 			& b&=\sqrt{AC}\; , 			& c&=\sqrt{AB}\; ,  & & \nonumber \\
x&=X\sqrt{BC}\; , 			& u&=U\sqrt{AC}\; , 			& v&=V\sqrt{AB}\; , & d&=\sqrt{EF}\; , \eqlabel{gencoordtrafo} \\ 
\alpha&=\mathcal{A}\sqrt{BC}\; , 	& y&=Y\sqrt{AC}\; , 			& w&=W\sqrt{AB}\; , & e&=\sqrt{DF}\; , \nonumber \\ 
\beta&=\mathcal{B}\sqrt{BC}\; , 	& \gamma&=\mathcal{C}\sqrt{AC}\; , 	& z&=Z\sqrt{AB}\; , & f&=\sqrt{DE} \nonumber \; .
\end{align}
One can easily check that this is, in fact, a diffeomorphism of $\left(\mathbb{R}_{>0}\right)^6\times \mathbb{R}^9$. In terms of the new coordinates, the expression for the scalar curvature $S$ is given by
\begin{align}
S&=BC+AC+AB+EF+DF+DE \nonumber \\
&+BC(X^2+\mathcal{A}^2+\mathcal{B}^2)+AC(U^2+Y^2+\mathcal{C}^2)+AB(V^2+W^2+Z^2) \nonumber \\
&-\frac{1}{2}\big(A^2+B^2+C^2+D^2+E^2+F^2\\
&+(B^2+C^2)(X^2+\mathcal{A}^2+\mathcal{B}^2)+(A^2+C^2)(U^2+Y^2+\mathcal{C}^2)+(A^2+B^2)(V^2+W^2+Z^2) \nonumber \\
&+(A(YZ-\mathcal{C} W)-DX)^2+(B(W\mathcal{B}-Z\mathcal{A})-DU)^2+(C(\mathcal{A}\mathcal{C}-\mathcal{B}Y)-DV)^2 \nonumber \\
&+(A(\mathcal{C} V-UZ)-E\mathcal{A})^2+(B(ZX-V\mathcal{B})-EY)^2+(C(\mathcal{B}U-X\mathcal{C})-EW)^2 \nonumber \\
&+(A(UW-YV)-F\mathcal{B})^2+(B(V\mathcal{A}-WX)-F\mathcal{C})^2+(C(XY-\mathcal{A}U)-FZ)^2\big) \; . \nonumber
\end{align}
In contrast to the rational expression~\eqref*{Sgen}, this is a polynomial of degree 6.
The volume $V$ in the old and new coordinates is given by
\begin{equation}
 V = (abcdef)^{-1} = (ABCDEF)^{-1} \; ,
\end{equation}
respectively.

\section{Left-invariant Einstein metrics invariant under a finite subgroup of \texorpdfstring{$\mathrm{Ad} (G)$}{Ad(G)}}

In this section we analyze left-invariant Einstein metrics $g$ on $G=S^3 \times S^3$ invariant under a non-trivial finite subgroup $\G\subset\mathrm{Ad}(G)$. We begin by observing that either all non-trivial elements of $\G$ are of order $2$ or there exists an element $\s$ of order $k\ge 3$. 
Let us first  consider the latter case. 
\bp Let $g$ be a left-invariant and $\G$-invariant Einstein metric on $G$, where $\G \subset \mathrm{Ad}(G)$.
If $\G$ contains an element $\s$ of order $k\ge 3$ then $K$, as defined in~\eqref*{Kdef},  contains a $U(1)$ subgroup and, hence, 
$(G,g)$ is homothetic to $(G,g_{can})$ or $(G,g_{NK})$.
\ep 

\pf 
Since $\mathrm{Ad} (G)$ is compact there exists a one-parameter subgroup 
which contains $\s$. Every one-parameter subgroup of $\mathrm{Ad} (G)$ is contained in a maximal 
torus $T \cong S^1\times S^1$ and $T$ is a product $S^1_1\times S^1_2$ of circle subgroups $S^1_1, S^1_2$ of the first and second $\SO (3)$-factors of
$\mathrm{Ad} (G)=\SO(3) \times \SO (3)$, respectively. Notice that the $\mathrm{Ad} (G)$-module $\gg = \mathrm{Lie}\, G$ is a sum 
\be 
 \gg = \bR^3_1 \oplus \bR^3_2
\ee
of $2$ inequivalent irreducible three-dimensional submodules $\mathbb{R}^3_\alpha$, $\alpha =1,2$, where the first factor of $\mathrm{Ad} (G)$ acts trivially on $\bR^3_2$ 
and the second factor acts trivially on $\mathbb{R}^3_1$. As $T$-modules we can decompose $\RR^3_\a$ further as 
\be
 \mathbb{R}^3_\alpha = \bR^1_\alpha \oplus \bR^2_\alpha \; ,
\ee
where $\bR^1_\alpha$ is a trivial module and $\bR^2_\alpha$ is irreducible.  
It follows that $\s$ acts as a rotation (with respect to the canonical scalar product) of order $k_\alpha\ge 1$ on $\bR^2_\alpha$, where 
$k_\alpha$ divides $k$ and at least one of the $k_\alpha$ is $\ge 3$, say $k_2 \ge 3$.   

If $k_1\neq k_2$ then the $\langle \s \rangle$-module $V^4:=\bR^3_1\oplus \bR^1_2=\bR^1_1\oplus \bR^2_1\oplus \bR^1_2$ does not contain any irreducible submodule 
equivalent to $\bR^2_2$. This implies that the submodules $V^4$ and $\bR^2_2\subset \gg$ are perpendicular for every $\langle \s \rangle$-invariant scalar product
on $\gg$.  Since $\s$ acts as a rotation of order $k_2\ge 3$ on $\bR^2_2$ it follows that the subgroup
$\SO (\bR^2_2) \subset \SO (\bR^3_2)=\{ e\} \times \SO(3) \subset \mathrm{Ad} (G)$ preserves every $\langle \s \rangle$-invariant scalar product
on $\gg$. Then the claim follows from \thmref{NikonorovThm}.  

If $k_1=k_2$, then $k_1=k_2=k$. In this case $\bR^2_1 \oplus \bR^2_2$ is the sum of $2$ equivalent 
irreducible $\langle \s \rangle$-modules and every $\langle \s \rangle$-invariant scalar product
on $\gg$ is invariant under the diagonally embedded subgroup $S^1\subset \SO (\bR^2_1) \oplus \SO(\bR^2_2)$
that contains $\s$. Thus, again, the claim follows from \thmref{NikonorovThm}.  \epf

For the remaining case we have the following result. 
\bp If all non-trivial elements of $\G$ are of order $2$, then   
$\G\cong \bZ_2^\ell$, where
$1\le \ell\le 4$. If $\ell \ge 3$,  then $\G$ contains an element $\s$ with $\tr \s = 2$.
\ep 

\pf 
Notice first that  $\G$ preserves the decomposition $\gg = \bR^3_1 \oplus \bR^3_2$. Moreover for given $\alpha \in \{ 1,2\}$, every non-trivial element $\s\in \G$ acts either 
trivially on $\bR^3_\alpha$ or $\bR^3_\alpha = \bR^1_\alpha \oplus \bR^2_\alpha$ is the sum of a trivial 
$\langle \s \rangle$-module and a non-trivial isotypical $\langle \s \rangle$-module (since $\s$ preserves the orientation
of $\bR^3_\alpha$), on which $\s$ acts as multiplication by $-1$. More precisely, either $\tr \s = -2$ or $\tr \s =2$ depending on whether
the $\langle \s \rangle$-modules  $\bR^3_1$ and $\bR^3_2$ are equivalent or not. 
It follows that $\ell \le 4$. If $\ell =4$ then $\gg$ splits as a sum of $6$ pairwise inequivalent one-dimensional $\G$-submodules. 
The last statement of the proposition is proven by simple combinatorics. \epf

The cases $\tr \s = 2$ and $\tr \s =-2$ will be treated separately. 
\bp  Let $g$ be a left-invariant and $\G$-invariant Einstein metric on $G$.  If $\G$ contains an involution $\s$ of trace $2$, then $g=g_{can}$. (By the previous proposition, this covers the case $\G\cong \bZ_2^\ell$, where
$\ell\ge 3$.)
\ep 

\pf 
We can assume that the $\langle \s \rangle$-module $\bR^3_1=\bR^1_1\oplus \bR^2_1$ is a sum of a trivial 
one-dimensional module and a nontrivial isotypical module, whereas $\bR^3_2$ is trivial. 
We show that the only left-invariant Einstein metric with normalized volume invariant under such an element is the standard metric. 
For every such metric $g$ there exists a $g$-orthonormal basis $(X_1,X_2,X_3,Y_1,Y_2,Y_3)$ such that 
$X_1\in \bR^1_1$, $X_2,X_3\in \bR^2_1$, and $Y_1, Y_2, Y_3 \in \bR^3_2$. Therefore it can be brought to the following form 
\be
 (X_1,X_2,X_3,Y_1,Y_2,Y_3) = (aE_1,bE_2,cE_3,xE_1+dF_1,\alpha E_1 + eF_2,\beta E_1 + fF_3) \; ,
\ee
with the same notation as introduced in \secref{gen}. Comparing with \eqref{genbasischange}, we learn that the equation above corresponds to the case where $y=z=u=v=w=\g=0$.

The scalar curvature~\eqref*{Sgen} thus simplifies to: 
\begin{eqnarray} 
 -2S&=&-2(a^2+b^2+c^2+d^2+e^2+f^2+x^2+\a^2+\b^2) \nonumber\\ 
 &&+ a^2b^2c^{-2}+b^2c^2a^{-2}+c^2a^2b^{-2} +d^2e^2f^{-2}+e^2f^2d^{-2}+f^2d^2e^{-2} \\
 &&+(b^2c^{-2}+c^2b^{-2})(x^2+\a^2+\b^2) + a^{-2}(d^2e^2f^{-2}\b^2 + d^2f^2e^{-2}\a^2 + e^2f^2d^{-2}x^2). \nonumber
\end{eqnarray}
For this $S$ we need to solve the variational problem~\eqref*{var3}. 
We first compute 
\be
 -2\frac{\partial S}{\partial x} = 2 x \left(\frac{ a^2 d^2 \left(b^2-c^2\right)^2+b^2 c^2 e^2 f^2}{a^2 b^2 c^2 d^2}\right) \; .
\ee
Notice that, since $a, \ldots, f$ are positive, the expression in parenthesis is positive. Therefore $\frac{\partial S}{\partial x} = \mu 
\frac{\partial V}{\partial x}=0$ implies that $x=0$. The same argument shows that $\a = \b =0$. 
Now the equation  $\nabla S = \mu \nabla V$ is equivalent to
\be\eqlabel{trace_plus_two_eqsys}
 a \frac{\p S}{\p a} = \cdots = f\frac{\p S}{\p f} = -\mu \; .
\ee
Together with the constraint equation $V=1$, this yields a system of 7 polynomial equations in the 7 unknowns $(a, \ldots, f, \mu)$ of degree at most 11,
\begin{align}
 0 &= a b c d e f-1 \; , \nonumber \\
 0 &= a b c \mu +a^4 b^4 d e f-a^4 c^4 d e f -b^4 c^4 d e f +2 a^2 b^2 c^4 d e f \; , \nonumber \\
 0 &= a b c \mu -a^4 b^4 d e f+a^4 c^4 d e f -b^4 c^4 d e f +2 a^2 b^4 c^2 d e f \; , \nonumber \\
 0 &= a b c \mu -a^4 b^4 d e f-a^4 c^4 d e f +b^4 c^4 d e f +2 a^4 b^2 c^2 d e f \; , \\
 0 &= d e f \mu +a b c d^4 e^4-a b c d^4 f^4-a b c e^4 f^4 +2 a b c d^2 e^2 f^4 \; , \nonumber \\
 0 &= d e f \mu -a b c d^4 e^4+a b c d^4 f^4-a b c e^4 f^4 +2 a b c d^2 e^4 f^2 \; , \nonumber \\
 0 &= d e f \mu -a b c d^4 e^4-a b c d^4 f^4+a b c e^4 f^4 +2 a b c d^4 e^2 f^2 \; . \nonumber
\end{align}
Manually solving this complicated system of coupled polynomial equations is unfeasible. Fortunately however it is well-suited for a computer-based Gr\"obner basis computation. (For an introductory text on the theory of Gr\"obner bases, see, for example,~\cite{CLO2015}.) As a result of such a Gr\"obner basis computation\footnote{This is the first
such computation in this paper, which is simple enough to be performed using any state of the art computer algebra software such as Mathematica, without additional hardware requirements. For the later calculations we will need more specific hardware and software, as described below.}, we find~\cite{GBurl}
\be
 a=b=c=d=e=f=-\mu=1
\ee
as the only solution with $a,\ldots,f\in \bR_{>0}$. This proves that $g=g_{can}$ if $\G$ contains an element of trace $2$. 
\epf 

It remains to treat the case when $1\le \ell\le 2$ and $\tr \s = -2$ for all non-trivial elements $\s\in \G$.  
In the following subsection, we first consider the case $\ell =2$.

\subsection{The case \texorpdfstring{$\G \cong \ZZ_2 \times \ZZ_2$}{Gamma = Z2xZ2}}\seclabel{Z2timesZ2}

When $\ell =2$ and all non-trivial elements of $\G\cong \ZZ_2 \times \ZZ_2$ are of trace $-2$, the $\G$ modules $\bR^3_1$, $\bR^3_2$ are equivalent and each of them splits as a sum
of three pairwise inequivalent one-dimensional submodules. This implies that there exists a 
$g$-orthonormal basis of the form 
\be
 (X_1,X_2,X_3,Y_1,Y_2,Y_3) = (aE_1,bE_2,cE_3,xE_1+dF_1,y E_2 + eF_2,z E_3 + fF_3) \; ,
\ee
where $a,\ldots,f\in \bR_{>0}$, $x, y, z \in \bR$, and $V=(abcdef)^{-1}=1$. 
Comparing with \eqref{genbasischange}, we learn that this corresponds to the case where $u=v=w=\a=\b=\g=0$. In this case the scalar curvature~\eqref*{Sgen} becomes
\begin{eqnarray} -2S&=&-2(a^2+b^2+c^2+d^2+e^2+f^2+x^2+y^2+z^2) \nonumber\\
&&+ a^2b^2c^{-2}+b^2c^2a^{-2}+c^2a^2b^{-2} +d^2e^2f^{-2}+e^2f^2d^{-2}+f^2d^2e^{-2} \nonumber\\
&&+y^2(a^2c^{-2}+c^2a^{-2}) + z^2 (a^2b^{-2}+b^2a^{-2})+x^2(b^2c^{-2}+c^2b^{-2}) \nonumber\\
&&+a^{-2}(yz -efd^{-1}x)^2 + b^{-2}(xz-dfe^{-1}y)^2 + c^{-2}(xy-def^{-1}z)^2. 
 \end{eqnarray}
For this $S$ we need to solve the variational problem~\eqref*{var3}, which we will achieve by again resorting to a computer-based Gr\"obner basis computation.

Before doing so, it is beneficial, in order to minimize the running time and complexity of the Gr\"obner basis computation, to utilize the coordinate transformation introduced in~\eqref*{gencoordtrafo}. With the simplification $u=v=w=\a=\b=\g=0$, the transformation only acts on the remaining coordinates $(a, \ldots, f, x, y, z) \to (A, \ldots, F, X, Y, Z)$ of $\left(\mathbb{R}_{>0}\right)^6\times \mathbb{R}^3$ and reads as follows
\begin{align}
 a&=\sqrt{B C} \; , &b&=\sqrt{A C} \; , &c&=\sqrt{A B} \; , \nonumber \\
 d&=\sqrt{F E} \; , &e&=\sqrt{D F} \; , &f&=\sqrt{D E} \; , \eqlabel{Z2xZ2coordtrafo} \\
 x&=X \sqrt{B C} \; , &y&=Y \sqrt{A C} \; , &z&=Z \sqrt{A B} \; . \nonumber
\end{align}
In terms of the new coordinates, the expression for the scalar curvature $S$ can be read off from
\begin{align}
 -2 S = &\quad A^2 + B^2 + C^2 + D^2 + E^2 + F^2 \nonumber \\ 
 & - 2 A B (1 + Z^2) - 2 A C (1 + Y^2) - 2 B C (1 + X^2) - 2 D E - 2 D F - 2 E F \\ 
 & + A^2 (Y^2 + Z^2) + B^2 (X^2 + Z^2) + C^2 (X^2 + Y^2) + D^2 X^2 + E^2 Y^2 + F^2 Z^2 \nonumber \\
 & - 2 (A D + B E + C F) X Y Z + A^2 Y^2 Z^2 + B^2 X^2 Z^2 + C^2 X^2 Y^2 \; . \nonumber
\end{align}
The variational problem~\eqref*{var3} leads altogether to ten polynomial equations of degree six in the ten unknowns $(A, \ldots, F, X, Y, Z, \mu)$,
\begin{align}
 0 &= A B C D E F-1, \nonumber \\
 0 &= B C D E F \mu -A Y^2 Z^2+D X Y Z-A Y^2-A Z^2 +B Z^2+C Y^2-A+B+C , \nonumber \\
 0 &= A C D E F \mu -B X^2 Z^2+E X Y Z-BZ^2-BX^2+C X^2+A Z^2+A-B+C  , \nonumber \\
 0 &= A B D E F \mu -C X^2 Y^2+F X Y Z-C X^2-C Y^2+A Y^2+B X^2+A+B-C  , \nonumber \\
 0 &= A B C E F \mu +A X Y Z-D X^2-D+E+F  , \nonumber \\
 0 &= A B C D F \mu +B X Y Z-E Y^2+D-E+F  , \eqlabel{Z2xZ2eqsys} \\
 0 &= A B C D E \mu +C X Y Z-F Z^2+D+E-F  , \nonumber \\
 0 &=-B^2 X Z^2-C^2 X Y^2-B^2 X-C^2 X+ A D Y Z+B E Y Z+C F Y Z+2 B C X-D^2 X  , \nonumber \\
 0 &= -C^2 X^2 Y-A^2 Y Z^2-C^2 Y-A^2 Y+A D X Z+B E X Z+C F X Z+2 A C Y-E^2 Y  , \nonumber \\
  0 &= -A^2 Y^2 Z-B^2 X^2 Z-A^2 Z-B^2 Z+A D X Y+B E X Y+C F X Y+2 A B Z-F^2 Z  . \nonumber
\end{align}
The polynomials on the right-hand sides form the input set for our Gr\"obner basis computation. We used the computer algebra system Magma~\cite{Magma1,Magma2} to compute\footnote{The computation was performed on a compute-server with 24 Intel Xeon E5-2643 3.40 GHz processors and 512 GB of RAM. The computational complexity is sensitive to the order of variables. We chose the following order of variables: $(A,B,C,D,E,F,X,Y,Z,\mu)$. The computation then took 16.5 minutes to run and consumed about 1.8 GB of RAM.\label{fref_hardware}} a Gr\"obner basis with lexicographic monomial ordering. 

The resulting Gr\"obner basis contains 55 polynomials with on average 78.7 terms per polynomial~\cite{GBurl}. The numerical coefficients range up to order $10^{12}$. Despite this apparent complexity, it is straightforward to find the vanishing locus of these polynomials owing to the elimination property of the lexicographic monomial ordering (see, for example,~\cite{CLO2015}). In terms of the original set of variables $(a, \ldots, f, x, y, z, \mu)$ we find \emph{a priori} 7 types of real solutions, as summarized in the following table.
\begin{center}\label{TabelleZ2xZ2}\def\arraystretch{1.6}
\begin{tabular}[h]{c|c|c|c|c|c|c|c|c|c|c|c}
 counter & $a$ & $b$ & $c$ & $d$ & $e$ & $f$ & $x$ & $y$ & $z$ & $\mu$ & $S$ \\ \hline
 $(1)$ & $1$ & $1$ & $1$ & $1$ & $1$ & $1$ & $0$ & $0$ & $0$ & $-1$ & $3$ \\
 $(2)$ & $1$ & $1$ & $1$ & $1$ & $1$ & $1$ & $\pm 1$ & $\pm 1$ & $1$ & $-1$ & $3$ \\
 $(3)$ & $1$ & $1$ & $1$ & $1$ & $1$ & $1$ & $\pm 1$ & $\mp 1$ & $-1$ & $-1$ & $3$ \\
 $(4)$ & $\frac{1}{\sqrt{2}}$ & $\frac{1}{\sqrt{2}}$ & $\frac{1}{\sqrt{2}}$ & $\sqrt{2}$ & $\sqrt{2}$ & $\sqrt{2}$ & $\pm\frac{1}{\sqrt{2}}$ & $\pm\frac{1}{\sqrt{2}}$ & $\frac{1}{\sqrt{2}}$ & $-1$ & $3$ \\
 $(5)$ & $\frac{1}{\sqrt{2}}$ & $\frac{1}{\sqrt{2}}$ & $\frac{1}{\sqrt{2}}$ & $\sqrt{2}$ & $\sqrt{2}$ & $\sqrt{2}$ & $\pm\frac{1}{\sqrt{2}}$ & $\mp\frac{1}{\sqrt{2}}$ & $-\frac{1}{\sqrt{2}}$ & $-1$ & $3$ \\
 $(6)$ & $\frac{\sqrt[4]{3}}{\sqrt{2}}$ & $\frac{\sqrt[4]{3}}{\sqrt{2}}$ & $\frac{\sqrt[4]{3}}{\sqrt{2}}$ & $\frac{\sqrt{2}}{\sqrt[4]{3}}$ & $\frac{\sqrt{2}}{\sqrt[4]{3}}$ & $\frac{\sqrt{2}}{\sqrt[4]{3}}$ & $\pm\frac{1}{\sqrt{2} \sqrt[4]{3}}$ & $\pm\frac{1}{\sqrt{2} \sqrt[4]{3}}$ & $\frac{1}{\sqrt{2} \sqrt[4]{3}}$ & $-\frac{5}{3 \sqrt{3}}$ & $\frac{5}{\sqrt{3}}$ \\ 
 $(7)$ & $\frac{\sqrt[4]{3}}{\sqrt{2}}$ & $\frac{\sqrt[4]{3}}{\sqrt{2}}$ & $\frac{\sqrt[4]{3}}{\sqrt{2}}$ & $\frac{\sqrt{2}}{\sqrt[4]{3}}$ & $\frac{\sqrt{2}}{\sqrt[4]{3}}$ & $\frac{\sqrt{2}}{\sqrt[4]{3}}$ & $\pm\frac{1}{\sqrt{2} \sqrt[4]{3}}$ & $\mp\frac{1}{\sqrt{2} \sqrt[4]{3}}$ & $-\frac{1}{\sqrt{2} \sqrt[4]{3}}$ & $-\frac{5}{3 \sqrt{3}}$ & $\frac{5}{\sqrt{3}}$ \\ 
\end{tabular}
\end{center}
Here, the first column represents a counter to distinguish the solutions, the last column contains the value of the scalar curvature $S$ at the respective solution (note that $S=6r$, with $r$ as defined in~\cite{NikRod}), and the signs in rows $2-7$ for $x$ and $y$ are correlated.

Note that the different choices of signs for the variables $x, y, z$ can be absorbed in the initial choice of the basis $(\mathbf{E},\mathbf{F})$, see 
above \eqref{genbasischange}. This reduces the above list to the four cases $(1)$, $(2)$, $(4)$, and $(6)$, with all the 
variables $x, y, z$ non-negative. 

We compare these solutions to the results already obtained in~\cite{NikRod} (in particular metrics $(1)$-$(4)$ in the proof of Lemma 2 on page 375). 
After adjusting notation, our solutions $(1)$, $(2)$, $(4)$, and $(6)$ correspond to the metrics $(1)$, $(2)$, $(3)$, and $(4)$ in~\cite{NikRod}, respectively. %
Our solution $(1)$ is the standard metric, $(6)$ is the nearly K\"ahler metric and $(2)$ and $(4)$ are isometric to 
the standard metric. The three metrics $(1)$, $(2)$, and $(4)$ correspond to the three possible decompositions 
of the manifold $S^3 \times S^3$ as a Riemannian product of two three-dimensional Lie subgroups (the two $S^3$-factors and 
the diagonal). 

We end this subsection by noting that altogether this completes the proof of \thmref{Z2xZ2}.

\subsection{The case \texorpdfstring{$\G \cong \ZZ_2$}{Gamma = Z2}}\seclabel{Z2}

In this subsection we consider the final remaining case, namely $\ell = 1$, that is $\G \cong \ZZ_2$, with the non-trivial element $\s\in\G$ satisfying $\tr \s = -2$. A qualitative novelty arises for this case, as will be explained below.

Fixing $\G \cong \ZZ_2$, with the non-trivial element $\s\in\G$ satisfying $\tr \s = -2$, implies that there exists a $g$-orthonormal basis of the form 
\be
 (X_1,X_2,X_3,Y_1,Y_2,Y_3) = (aE_1,bE_2,cE_3,xE_1+dF_1,y E_2 + w E_3 + eF_2,\g E_2 + z E_3 + fF_3) \; ,
\ee
where $a,\ldots,f\in \bR_{>0}$, $x, y, z, w, \g \in \bR$, and $V=(abcdef)^{-1}=1$. 
Comparing with \eqref{genbasischange}, we learn that this corresponds to the case where $u=v=\a=\b=0$ and the scalar curvature~\eqref*{Sgen} becomes
\begin{align}
 S &= a^2 +b^2 +c^2 + d^2 + e^2 + f^2 + x^2 + y^2 + z^2 +w^2 + \g^2 \nonumber \\
  &\quad - \frac12 \left\{ \vphantom{\left[ \left(\frac{de}{f}\right)^2 \right]} a^2 b^2 c^{-2} + b^2 c^2 a^{-2} + c^2 a^2 b^{-2} + d^2 e^2 f^{-2} + e^2 f^2 d^{-2} + f^2 d^2 e^{-2} \right. \nonumber \\ 
  &\quad + \left( \frac{a^2}{c^2} + \frac{c^2}{a^2} \right) (y^2 + \g^2) + \left( \frac{a^2}{b^2} + \frac{b^2}{a^2} \right) (w^2 + z^2) + \left( \frac{b^2}{c^2} + \frac{c^2}{b^2} \right) x^2 
  + a^{-2} \left( yz - w\g - \frac{ef}{d} x \right)^2
  \nonumber \\
  &\quad+ b^{-2} \left[ \left( xw + \frac{de}{f} \g \right)^2 + \left( xz - \frac{df}{e} y \right)^2 \right]
  + c^{-2} \left. \left[ \left( xy - \frac{de}{f} z \right)^2 + \left( x\g + \frac{df}{e} w \right)^2 \right] \right\} \; .
\end{align}
Next, we again employ the coordinate transformation~\eqref*{gencoordtrafo} in order to facilitate the upcoming Gr\"obner basis computation.
With the simplification $u=v=\a=\b=0$, the transformation only acts on the remaining coordinates $(a, \ldots, f, x, y, z, w, \g) \to (A, \ldots, F, X, Y, Z, W, \mathcal{C})$ of $\left(\mathbb{R}_{>0}\right)^6\times \mathbb{R}^5$,
\begin{align}
 a&=\sqrt{B C} \; , &b&=\sqrt{A C} \; , &c&=\sqrt{A B} \; , \nonumber \\
 d&=\sqrt{F E} \; , &e&=\sqrt{D F} \; , &f&=\sqrt{D E} \; , \eqlabel{Z2coordtrafo} \\
 x&=X \sqrt{B C} \; , &y&=Y \sqrt{A C} \; , &z&=Z \sqrt{A B} \; , \nonumber \\
 w&=W \sqrt{A B} \; , &\g&=\mathcal{C} \sqrt{A C}\; . && \nonumber
\end{align}
In terms of the new coordinates, the scalar curvature $S$ is given by
\begin{align}
 S &=  -\, \frac{A^2}{2} - \frac{B^2}{2} - \frac{C^2}{2} - \frac{D^2}{2} - \frac{E^2}{2} - \frac{F^2}{2}  + A B   +  A C + B C  + D F + D E + F E \nonumber \\
     &\quad - \frac{A^2 Y^2}{2} - \frac{A^2 Z^2}{2} - \frac{A^2 W^2}{2} - \frac{A^2 \mathcal{C} ^2}{2} - \frac{B^2 X^2}{2}  - \frac{B^2 Z^2}{2} - \frac{B^2 W^2}{2} \nonumber \\
     &\quad - \frac{C^2 X^2}{2}  - \frac{C^2 Y^2}{2} - \frac{C^2 \mathcal{C} ^2}{2} - \frac{D^2 X^2}{2} - \frac{E^2 Y^2}{2} - \frac{E^2 W^2}{2}  - \frac{F^2 Z^2}{2}  - \frac{F^2 \mathcal{C} ^2}{2} \nonumber \\
     &\quad + A B Z^2 + A B W^2 + A C \mathcal{C} ^2  + B C X^2  + A C Y^2 \nonumber \\
     &\quad + A D X Y Z + C F X Y Z + B E X Y Z  - A D W X \mathcal{C}  - B F W X \mathcal{C}  - C E W X \mathcal{C}  + A^2 W Y Z \mathcal{C} \nonumber \\
     &\quad - \frac{1}{2} B^2 W^2 X^2 - \frac{1}{2} B^2 X^2 Z^2 - \frac{1}{2} A^2 Y^2 Z^2  - \frac{1}{2} C^2 X^2 Y^2  - \frac{1}{2} A^2 W^2 \mathcal{C} ^2 - \frac{1}{2} C^2 X^2 \mathcal{C} ^2
\end{align}
The variational problem~\eqref*{var3} leads altogether to 12 polynomial equations of degree 6 in the 12 unknowns $(A, \ldots, F, X, Y, Z, W, \mathcal{C}, \mu)$:
\begin{align}
 0 &= A B C D E F - 1, \nonumber \\
 0 &=-D+E+F + A B C E F \mu-D X^2+A X Y Z-A W X \mathcal{C} , \nonumber \\
 0 &=D-E+F + A B C D F \mu-E W^2-E Y^2+B X Y Z-C W X \mathcal{C} , \nonumber \\
 0 &=A-B+C + A C D E F \mu+A W^2-B W^2-B X^2+C X^2-B W^2 X^2+E X Y Z \nonumber \\
    &+A Z^2-B Z^2-B X^2 Z^2-F W X \mathcal{C} , \nonumber \\
 0 &=D+E-F + A B C D E \mu+C X Y Z-F Z^2-B W X \mathcal{C} -F \mathcal{C} ^2, \nonumber \\    
 0 &=-A+B+C + B C D E F \mu-A W^2+B W^2-A Y^2+C Y^2+D X Y Z-A Z^2+B Z^2 \nonumber \\
    &-A Y^2 Z^2-D W X \mathcal{C} +2 A W Y Z \mathcal{C} -A \mathcal{C} ^2+C \mathcal{C} ^2-A W^2 \mathcal{C} ^2, \nonumber \\
 0 &=A+B-C + A B D E F \mu+B X^2-C X^2+A Y^2-C Y^2-C X^2 Y^2+F X Y Z-E W X \mathcal{C} \nonumber \\
    &+A \mathcal{C} ^2-C \mathcal{C} ^2-C X^2 \mathcal{C} ^2 \nonumber \\
 0 &=-A D W X-C E W X-B F W X+A^2 W Y Z-A^2 \mathcal{C} +2 A C \mathcal{C} -C^2 \mathcal{C} -F^2 \mathcal{C} \nonumber \\
    &-A^2 W^2 \mathcal{C} -C^2 X^2 \mathcal{C} , \nonumber \\
 0 &=A D X Y+B E X Y+C F X Y-A^2 Z+2 A B Z-B^2 Z-F^2 Z-B^2 X^2 Z \nonumber \\
    &-A^2 Y^2 Z+A^2 W Y \mathcal{C} , \nonumber \\
 0 &=-A^2 Y+2 A C Y-C^2 Y-E^2 Y-C^2 X^2 Y+A D X Z+B E X Z+C F X Z \eqlabel{Z2eqsys} \\
    &-A^2 Y Z^2+A^2 W Z \mathcal{C} , \nonumber \\
 0 &=-A^2 W+2 A B W-B^2 W-E^2 W-B^2 W X^2-A D X \mathcal{C} -C E X \mathcal{C} -B F X \mathcal{C} \nonumber \\
    &+A^2 Y Z \mathcal{C} -A^2 W \mathcal{C} ^2, \nonumber \\
 0 &=-B^2 X+2 B C X-C^2 X-D^2 X-B^2 W^2 X-C^2 X Y^2+A D Y Z+B E Y Z+C F Y Z \nonumber \\
    &-B^2 X Z^2-A D W \mathcal{C} -C E W \mathcal{C} -B F W \mathcal{C} -C^2 X \mathcal{C} ^2 \nonumber \; .
\end{align}
The polynomials on the right-hand sides form the input set for our Gr\"obner basis computation. Unfortunately, computing a Gr\"obner basis with lexicographic monomial ordering, and consequently solving the system, is apparently out of reach with current technology.

However, it is possible to compute a Gr\"obner basis with graded reverse lexicographic (or \emph{grevlex}, for short) monomial ordering, instead. This took about 29 days to run\footnote{See footnote~\ref{fref_hardware} for a description of the hardware used to perform the computation. The order of variables was in this case chosen to be $(\mu, D, F, E, C, B, Z, Y, \mathcal{C}, W, A, X)$.} and consumed about 78 GB of RAM. The generated output has a size of 106 GB in a human-readable format. It consists of 50472 polynomials with on average 593 terms per polynomial~\cite{GBurl}. The numerical coefficients range up to order $10^{10}$. 
Since the grevlex Gr\"obner basis lacks the elimination property, it is not helpful for solving the system, but can be used to examine general properties of the solution set. In particular, by applying the Finiteness Theorem of~\cite[p. 251, \S 5.3, Theorem 6]{CLO2015} one learns whether or not the solution set is finite (over the complex numbers). From the Finiteness Theorem we conclude that the system~\eqref*{Z2eqsys} has a continuous family of \emph{complex} solutions.

We remark that
this is a qualitative novelty compared to the other cases considered in this paper. Indeed, regarded as complex varieties,~\eqref*{trace_plus_two_eqsys,Z2xZ2eqsys} are zero-dimensional, whereas the dimension of the complex variety defined by~\eqref*{Z2eqsys} is larger than zero. 
This observation has consequences for the Gr\"obner basis computation, since more efficient algorithms are available for the case of zero-dimensional varieties. This technicality at least partly explains why we have not been able to compute a Gr\"obner basis with lexicographic monomial ordering for the system~\eqref*{Z2eqsys}.

The complexity, and hence running time, of Gr\"obner basis computations typically scales rather badly (that is, doubly exponentially) in terms of the size of the input, which is in turn related to the number of variables, number of polynomials, and degrees of the polynomials (see, for example,~\cite[\S 21.7]{vzGathen_Gerhard_2013}, and references therein, for a brief summary of the current status on the complexity of Gr\"obner basis computations). We may hope to be able to perform the desired computation of the Gr\"obner basis with lexicographic monomial ordering if we consider restrictions of the polynomial system~\eqref*{Z2eqsys}.

This is indeed the case if we fix, for example, the value of the Lagrange multiplier $\mu$. Two distinguished cases are $\mu=-1$ and $\mu=-5/(3 \sqrt{3})$, which correspond to the known solutions $(G,g_{can})$ and $(G,g_{NK})$ found in \thmref{NikonorovThm,Z2xZ2}.

In the first case, we add the polynomial $\mu + 1$ to the input set given by the right hand sides of~\eqref*{Z2eqsys} and compute the Gr\"obner basis with lexicographic monomial ordering for the variable ordering $(\mu,F,E,D,C,B,A,Z,W,\mathcal{C},Y,X)$. The computation takes about 76 minutes to run and consumed about 3.4 GB of RAM. The resulting Gr\"obner basis has a size of 551 bytes and consists of 16 polynomials~\cite{GBurl}.
In terms of the original set of variables $(a, \ldots, f, x, y, z, w, \g)$ we find \emph{a priori} 5 types of real solutions, as summarized in the following table.
\begin{center}\def\arraystretch{1.6}
\begin{tabular}[h]{c|c|c|c|c|c|c|c|c|c|c}
 $a$ & $b$ & $c$ & $d$ & $e$ & $f$ & $x$ & $y$ & $z$ & $w$ & $\g$ \\ \hline
 $1$ & $1$ & $1$ & $1$ & $1$ & $1$ & $0$ & $0$ & $0$ & $0$ & $0$ \\
 $1$ & $1$ & $1$ & $1$ & $1$ & $1$ & $-1$ & $t$ & $-t$ & $\pm\sqrt{1 - t^2}$ & $\pm\sqrt{1 - t^2}$ \\
 $1$ & $1$ & $1$ & $1$ & $1$ & $1$ & $1$ & $t$ & $t$ & $\mp\sqrt{1 - t^2}$ & $\pm\sqrt{1 - t^2}$ \\
 $\frac{1}{\sqrt{2}}$ & $\frac{1}{\sqrt{2}}$ & $\frac{1}{\sqrt{2}}$ & $\sqrt{2}$ & $\sqrt{2}$ & $\sqrt{2}$ & $-\frac{1}{\sqrt{2}}$ & $\frac{t}{\sqrt{2}}$ & $-\frac{t}{\sqrt{2}}$ & $\pm\frac{\sqrt{1-t^2}}{\sqrt{2}}$ & $\pm\frac{\sqrt{1-t^2}}{\sqrt{2}}$ \\
 $\frac{1}{\sqrt{2}}$ & $\frac{1}{\sqrt{2}}$ & $\frac{1}{\sqrt{2}}$ & $\sqrt{2}$ & $\sqrt{2}$ & $\sqrt{2}$ & $\frac{1}{\sqrt{2}}$ & $\frac{t}{\sqrt{2}}$ & $\frac{t}{\sqrt{2}}$ & $\mp\frac{\sqrt{1-t^2}}{\sqrt{2}}$ & $\pm\frac{\sqrt{1-t^2}}{\sqrt{2}}$
\end{tabular}
\end{center}
Here, the quantity $t\in[-1,1]$ is a free parameter and the signs in rows $2-5$ for $w$ and $\g$ are correlated. The value of the scalar curvature $S$ is $3$ for all of the above solutions. It can be shown by a change of  the initial basis $(\mathbf{E},\mathbf{F})$ that all solutions are isometric to the standard metric $g_{can}$, irrespective of the value of the parameter $t$, see the remarks after the table on page \pageref{TabelleZ2xZ2}.

In the second case, we add the polynomial $\mu + 5/(3 \sqrt{3})$ to the input set given by the right hand sides of~\eqref*{Z2eqsys} and compute the Gr\"obner basis with lexicographic monomial ordering for the variable ordering $(\mu,F,E,D,C,B,A,Z,W,\mathcal{C},Y,X)$. The computation takes about 19 minutes to run and consumed about 1.3 GB of RAM. The resulting Gr\"obner basis has a size of 535 bytes and consists of 12 polynomials~\cite{GBurl}.
In terms of the original set of variables $(a, \ldots, f, x, y, z, w, \g)$ we find \emph{a priori} 2 types of real solutions, as summarized in the following table.
\begin{center}\def\arraystretch{1.6}
\begin{tabular}[h]{c|c|c|c|c|c|c|c|c|c|c}
 $a$ & $b$ & $c$ & $d$ & $e$ & $f$ & $x$ & $y$ & $z$ & $w$ & $\g$ \\ \hline
 $\frac{\sqrt[4]{3}}{\sqrt{2}}$ & $\frac{\sqrt[4]{3}}{\sqrt{2}}$ & $\frac{\sqrt[4]{3}}{\sqrt{2}}$ & $\frac{\sqrt{2}}{\sqrt[4]{3}}$ & $\frac{\sqrt{2}}{\sqrt[4]{3}}$ & $\frac{\sqrt{2}}{\sqrt[4]{3}}$ & $\frac{1}{\sqrt{2} \sqrt[4]{3}}$ & $\pm\frac{\sqrt{1-3 t^2}}{\sqrt{2} \sqrt[4]{3}}$ & $\pm\frac{\sqrt{1-3 t^2}}{\sqrt{2} \sqrt[4]{3}}$ & $-\frac{\sqrt[4]{3}}{\sqrt{2}}\, t$ & $\frac{\sqrt[4]{3}}{\sqrt{2}}\, t$ \\
 $\frac{\sqrt[4]{3}}{\sqrt{2}}$ & $\frac{\sqrt[4]{3}}{\sqrt{2}}$ & $\frac{\sqrt[4]{3}}{\sqrt{2}}$ & $\frac{\sqrt{2}}{\sqrt[4]{3}}$ & $\frac{\sqrt{2}}{\sqrt[4]{3}}$ & $\frac{\sqrt{2}}{\sqrt[4]{3}}$ & $-\frac{1}{\sqrt{2} \sqrt[4]{3}}$ & $\mp\frac{\sqrt{1-3 t^2}}{\sqrt{2} \sqrt[4]{3}}$ & $\pm\frac{\sqrt{1-3 t^2}}{\sqrt{2} \sqrt[4]{3}}$ & $\frac{\sqrt[4]{3}}{\sqrt{2}}\, t$ & $\frac{\sqrt[4]{3}}{\sqrt{2}}\, t$ 
\end{tabular}
\end{center}
Here, the quantity $t\in[-1/\sqrt{3},1/\sqrt{3}]$ is a free parameter and in both rows the signs for $y$ and $z$ are correlated. The value of the scalar curvature $S$ is $5/\sqrt{3}$ for all of the above solutions. It can be shown by a change of the initial basis $(\mathbf{E},\mathbf{F})$  that all solutions are isometric to the nearly K\"ahler metric $g_{NK}$, irrespective of the value of the parameter $t$.

This completes the proof of \propref{Z2lambdafixed}.
We conclude that if a new left-invariant Einstein metric of normalized volume on $S^3 \times S^3$ with additional orientation preserving $\bZ_2$-symmetry exists, its scalar curvature is different from that of the two known examples.
The above calculations indicate that for every specified value of the scalar curvature it should be possible to decide, using Gr\"obner basis methods, whether a left-invariant Einstein metric of normalized volume and given value of the  scalar curvature exists. In fact, we have applied this method to a small number of other values of the scalar curvature and, in each case, found that no (complex) solution to~\eqref*{Z2eqsys} exists, with the exception of the cases $\mu \in \{ 0, 1, 5/(3 \sqrt{3}), \pm 2/\sqrt{3} \}$, which yield complex (but not real) solutions~\cite{GBurl}.

\begin{remark} The existence of the one-parameter families of solutions is due to the ambiguity of the normal form of the metric $g$ in the cases where the matrix $A$ defined in \eqref{genbasischange} has multiple eigenvalues. Indeed, in these cases one can use the freedom in the choice of the initial basis $(\mathbf{E},\mathbf{F})$ to reduce the number of parameters in the off-diagonal square matrix $W$,  see \eqref{genbasischangematrix} for the notations. Therefore, it is not surprising that, in the absence of such an {\em a priori} reduction of the number of variables, the system \eqref*{Z2eqsys} admits one-parameter families of solutions, as shown in the two tables above. In fact, conjugating the matrix $A$ (encoding the solution) by a one-parameter group of rotations commuting with the diagonal part of $A$ produces a one-parameter family of isometric solutions. In this way, one can even obtain families depending on more than one parameter\footnote{For example, applying this observation to the metric $(2)$ from the table on page~\pageref{TabelleZ2xZ2}, we obtain a three-dimensional family of metrics, which are all isometric to the product metric. We thank Yuri Nikonorov for pointing this out to us.}, which are however not automatically in the considered normal form for $\bZ_2$-invariant metrics. Bringing these metrics to the normal form reduces the number of parameters.
It is an open question if the  algebraic subset  of  $\left(\mathbb{R}_{>0}\right)^6\times \mathbb{R}^6$ defined by \eqref*{Z2eqsys} can be decomposed into its intersections with  the two hyperplanes $\{\mu=-1\}$ and $\{\mu=-5/(3\sqrt 3)\}$ and some additional finite set (for which one can hope to determine all its points by computer algebra methods).\end{remark}

\resumetocwriting

\end{document}